\documentclass[10pt]{article}
\usepackage{amsfonts}
\usepackage{amsfonts}  
\title{Superlinear problems related to Fully Non Linear singular or degenerate operators in  punctured balls }

\author{  Fran\c{c}oise Demengel\\
  D\'epartement de Math\'ematiques,
CY Paris University
 }

\date{}

\usepackage{amsthm}
\usepackage{amsmath}
\usepackage{color}

\usepackage[colorlinks=true]{hyperref}

\catcode`@=11 \@addtoreset{equation}{section} \catcode`@=12

\newtheorem{theo}{Theorem}[section]
\newtheorem{prop}[theo]{Proposition}
\newtheorem{rema}[theo]{Remark}

\newtheorem{lemme}[theo]{Lemma}

\def\R{\mathbb  R}

\begin{document}
\maketitle

\begin{abstract}
We study Fully Nonlinear  uniformly elliptic equations with a higher differential term degenerate or singular, coupled with a zero order term singular at the point $0$, and an absorbing super-linear term of $p$ power type. 

 We will limitate to punctured balls and positive  radial solutions. In the continuity of two previous works in \cite {BDL1}, \cite{BDL2}, we will provide a complete classification depending on $p$ of the asymptotic behaviour of the solutions near the singular point $0$.
 \end{abstract}

 \section{Introduction}
 
  We are interested in the existence and behavior around the origin of   positive radial solutions of Fully Non Linear elliptic equations posed in punctured balls centered at the origin, namely equations of the form 
  \begin{equation}\label{eq1}
  | \nabla u |^\alpha F( D^2 u) + \mu {u^{1+\alpha} \over |x|^{2+\alpha}} = u^p. 
  \end{equation}
  
  Here  $\alpha >-1$, $p> 1+\alpha$, $\mu$ is a positive parameter to be precised later, and $F$ is one of the Pucci's operators, even if  some of the results here enclosed are  also valid for more general FNL operators. The definition of $| \nabla u |^\alpha F( D^2 u)$ is intended
 in the sense of \cite{BD1}, \cite{BD2}. Note that this type of operators is an extension of quasi linear operators as the $p$Laplacian operator with $p= 2+\alpha$ to a fully non linear framework.

 When $F$ is the Laplace operator and $\alpha = 0$  problem \eqref{eq1} has been extensively studied, also in the case $\mu=0$, see \cite{Ci, CirsChau, CDu1, CDu2, WDu}.  However, the techniques  employed rely on the variational structure of the problem, and they cannot  immediately be  extended to the fully nonlinear framework.
  Let us recall the definition of the extremal Pucci's operators 
\begin{equation}\label{M+} \mathcal{M}^+_{\lambda, \Lambda} (M)= \Lambda \sum_{\mu_i>0}\mu_i+\lambda\, \sum_{\mu_i<0}\mu_i=\sup_{A\in {\mathcal A}_{\lambda,\Lambda}} {\rm tr} (AM)\, ,\mathcal{M}^-_{\lambda, \Lambda} (M)= \lambda \sum_{\mu_i>0}\mu_i+\Lambda\, \sum_{\mu_i<0}\mu_i
=\inf_{A\in {\mathcal A}_{\lambda,\Lambda}} {\rm tr} (AM)\, ,
 \end{equation}
 where $\mu_1,\ldots,\mu_N$ are the eigenvalues of any  matrix $M$ belonging to the space $\mathcal{S}_N$ of the $N-$squared symmetric matrices.
In the sequel we will suppose $\lambda$ and $\Lambda$ fixe and denote $ \mathcal {M}^+( M) $  for simplicity  in place of $ \mathcal{M}^+_{\lambda, \Lambda} (M)$. 
When $u$ is a radial ${\cal C}^2$ function in $\R^N\setminus \{0\}$, using the abuse $u(x) = u(|x|) = u(r)$,   the eigenvalues of $D^2 u$ are $u^{\prime\prime}$ of multiplicity $1$ and ${ u^\prime \over r}$ of multiplicity $N-1$. As a consequence for example when $u$ is decreasing and convex ( which will almost always be the case here).
  $$ { \cal M}^+( D^2 u) = \Lambda u^{\prime \prime } + \lambda ( N-1) { u^\prime (r) \over r}.$$
  
  The purely eigenvalue problem associated with the operator $| \nabla u |^\alpha {\cal M}^\pm ( D^2 u)$ and the singular potential $r^{-2-\alpha}$ ie 
   $$\left\{ \begin{array}{lc}\label{eqM}
    | \nabla u |^\alpha { \cal M}^\pm ( D^2 u) + \lambda^\pm { u^{1+\alpha} \over r^{2+\alpha} }= 0 &{\rm in} \ B_1\setminus\{0\} \\
    u=0& {\rm on } \ \ \partial B_1.\end{array}\right.$$
     has been recently studied in \cite{D1}. The case $\alpha = 0$ is treated in \cite{BDL1} .  Let us recall the definitions of 
     $N^+$ and $N^-$ which are the $N$ dimension adapted to Pucci's operators, see for example  Cutri and Leoni \cite{CL}:
     
      $$N^+ = { \lambda \over \Lambda } ( N-1)+1, N^- = { \Lambda \over \lambda } ( N-1)+1.$$
We explicitely determined  in \cite{D1}
     $$\bar \lambda^+= \Lambda \left({N^+ -1\over 2+\alpha}\right)^{2+\alpha} (1+\alpha)^{1+\alpha}, $$
     in the case of ${ \cal M}^+$ while in the case of ${ \cal M}^-$ we exchange   $\lambda$ and $\Lambda$,   $N^+$ and $N^-$.

     In the sequel we will focus on the case of ${ \cal M}^+$, the changes to bring to the case of ${\cal M}^-$ being simple and left to the reader. For this reason we have denoted $\bar \lambda$ in place of $\bar \lambda^+$.  The equation we will then refer  to is  
       \begin{equation}\label{eqmain}
  | \nabla u |^\alpha {\cal M}^+( D^2 u) + \mu {u^{1+\alpha} \over |x|^{2+\alpha}} = u^p. 
  \end{equation}

      By exploiting the invariance of the equation with respect to the scaling
 $$
 u_\gamma (r)=\gamma ^{\frac{2+\alpha }{1+\alpha -p}} u \left( \frac{r}{\gamma }\right)\qquad \hbox{ for all } \gamma>0\, ,
 $$
 we will consider without loss of generality equation \eqref{eqmain} posed in 
 $ B_1\setminus \{0\}$, 
 and we will refer to smooth, positive, radial solutions  simply as solutions.

      In order to state our main results, let $\mu < \bar \lambda $ and let us introduce $\tau^\pm$ the two  positive zeros of  
      
       $$X^{2+\alpha} -(N^+-2) X^{1+\alpha} + {\mu\over \Lambda}.$$
        We have 
        $$\tau^-< \tau := (1+\alpha) ( N^+-2)(2+\alpha)^{-1}< \tau^+< N^+-2.$$
      We need also to  introduce the critical values 
        $p^\star = 1+\alpha + {2+\alpha \over \tau^+}$, $p^{\star\star} = 1+\alpha + {2+\alpha \over \tau^-}$. Note that 
        $1+\alpha < p^\star< p^{ \star\star}
$

 Our first result is an existence's result of radial positive solutions having the desired behavior near zero. 
 Let us define for simplicity the value 
 \begin{equation}\label{beta}
  \beta = {2+\alpha \over p-1-\alpha}.
  \end{equation}
 \begin{theo}\label{exi}
  If $\mu < \bar \lambda$ and $p> 1+\alpha$, 
  
  \begin{enumerate}
  \item If $p< p^\star$ the equation   \eqref{eqmain} has at least one solution  so that 
  $$r^{\beta }    u(r) \rightarrow K, $$
  where 
$$K^{p-1-\alpha} =\left( \Lambda \beta ^{\alpha+1} (\beta-( N^+-2))+ \mu\right)$$

  and for any $c^\pm >0$ at least solutions $u$ satisfying 
  $$r^{ \tau^\pm} u(r) \rightarrow c^\pm.$$
   \item 
   If $p^\star \leq p< p^{\star\star}$, 
   for any $c>0$,  there exists $u$ , solution of \eqref{eqmain} satisfying
    $$   u(r) \sim c r^{-\tau^-}.$$
    
    \item If $p = p^{\star \star}$ there exists at least one solution  of \eqref{eqmain} satisfying 
    $$r^{\tau^-} ( -\log r)^{ \tau^-\over 2+\alpha}     u(r) \rightarrow \tilde K, $$
     where 
     $$\tilde K = ( \Lambda (\tau^- )^{1+\alpha} ( \tau-\tau^-))^{\tau^-\over 2+\alpha}.$$
     
     \item If $p> p^{\star\star}$ there exists at least one solution of \eqref{eqmain}satisfying
    $$r^{\beta}    u(r) \rightarrow K.$$
    \end{enumerate}
    \end{theo}
    The second Theorem says that there is no other possible behaviors than those described in Theorem \ref{exi}
    
    \begin{theo}\label{behav}
  If $\mu < \bar \lambda$ and $p> 1+\alpha$, let $u$ be  a positive radial solution of \eqref{eqmain} in $B_1\setminus \{0\}$. Then 
  
  \begin{enumerate}
  \item If $p< p^\star$,  $u$ satisfies either 
    $$r^{\beta}    u(r) \rightarrow K, $$
 or there exist some positive constants $c^\pm$ so that 
  $$r^{ \tau^\pm} u(r) \rightarrow c^\pm.$$
   \item 
   if $p^\star \leq p< p^{\star\star}$, 
   $u$ satisfies  for some positive constant $c$ 
    $$  u(r) \sim c r^{-\tau^-}.$$
    
    \item If $p = p^{\star \star}$,  $u$ satisfies 
    $$r^{\tau^-} ( -\log r)^{ \tau^-\over 2+\alpha}     u(r) \rightarrow \tilde K.$$
         
     \item If $p> p^{\star\star}$     $$r^{\beta }    u(r) \rightarrow K.$$
    \end{enumerate}
    \end{theo}

      \bigskip

      Even if the Emdem Fowler transform ( see \cite{Har} for example)  does not permit to make explicit calculations of the solutions as in the case $\alpha = 0$, (see \cite{BDL2})  we will frequently  use the change of unknown 
      $$x(t) = r^{\beta} u(r), r = e^t.$$
      In particular the convergence $r^{\beta } u(r)\rightarrow 0$ means $x(t) \rightarrow 0$ at $-\infty$. 
      This change of functions reduces the initial problem for $u$ to     
     a two dimensional  differential system of  order one : 
    
     $$\left\{\begin{array}{lc} 
      x^\prime = \beta x + |Y|^{{1\over 1+\alpha}-1} Y& \ \\
      {Y^\prime \over 1+\alpha} = (\beta+1) Y + {1\over \Lambda} \left( -\lambda ( N-1) Y + x^p-\mu x^{1+\alpha} \right)^+ -{1\over \lambda}  \left(- \lambda ( N-1) Y + x^p-\mu x^{1+\alpha} \right)^- &
      \end{array}\right.$$

\bigskip

    The proofs of Theorem \ref{exi} and \ref{behav} will be obtained with the aid of some preliminary results : 
    
    In section 2 we give two comparison Theorems : The first one is valid not only for the Pucci's operators, nor only   for radial solutions.     
    In section 3 we prove the existence of solutions with the desired behavior. In section 4, after some preliminary observations on the behavior of non trivial solutions we prove that any solution satisfies the behavior described in the existence's Theorem.  An essential difference with the case $\alpha \neq 0$ treated in \cite{BDL2} is    that we will need  technical propositions  which   permit to " linearize" in some sort the  higher  order  part  of derivation of the  equation, in order to  provide   a complete  asymptotic behavior when $u(r) r^\beta$ goes to zero.

 In conclusion, let us finally mention that the results of the present paper can be regarded as a first step into the study   of fully nonlinear uniformly elliptic equations having   coefficients with isolated singularities. Many  questions remain open, ranging from the analysis of non radially symmetric solutions to the extension of the results to more general elliptic operators and zero order terms. In particular, the study of the limit case  $\mu=\bar \lambda$ is  postponed to a future work, as well as the case $\mu > \bar \lambda$. 
 
 \bigskip
 
  Note that a big amount of papers are interested in the  case where $u^p$ is replaced by $-u^p$. One good reference  in  the case of the Laplacian is \cite{Ci}, while we  consider  in \cite{BDL3} the Fully Non linear case  with $\alpha = 0$. The case of Fully Non linear operators,  treated in \cite{FQ} in the case  $\mu = 0$, gives an idea of the  problems encountered, even if  here the singularity at zero  brings  additional  difficulties.   One  difference with the "absorbing" case,   is the existence of periodic orbits  for some particular values of $p$, more precisely  in the case $p \in [p^\star, {4\over N^+-2}]$. We can expect to have  similar results  when $\alpha \neq 0$, but the highest non linearity of the  equation  brings new technical difficulties. This will be the object of a future work.

\section{Preliminary results}
We begin with the two comparison Theorems announced in the introduction : 
 
               \begin{theo}\label{comp1}
               
                Suppose that $u\geq 0$ and $v$   positive are respectively   viscosity sub and super-solutions of the equation 
                
                $$| \nabla u |^\alpha F( D^2 u) + \mu u^{1+\alpha}  r^{-2-\alpha } = u^p$$                 with $p > 1+\alpha$,                  and both are  continuous            
                 on the annulus $]r_1,  r_2[$ with $r_1>0$. Then if $u\leq v$ on $r_i$, $i = 1,2$ one has 
                 $$u\leq v$$ in all the annulus.
                 
                  \end{theo}
                  
                   \begin{proof}
                   
                    Let us suppose by contradiction that $\sup_{|r_1, r_2]} { u \over v} = \gamma^\prime >1$. Let us take $ \gamma \in ]1, \gamma^\prime[$ close enough to $\gamma^\prime$ so that 
                    
                     $$\mu (\left({\gamma^\prime \over \gamma}\right)^{1+\alpha}-1) \leq \inf_{ [r_1, r_2]}( r^{2+\alpha} v^{p-1-\alpha}) ( \gamma^{p-1-\alpha}-1).$$
                      The supremum of $u-\gamma v$ is positive and achieved inside the   open intervall $]r_1, r_2[$. 
                      Then using the fact that $u$ and $v$ are sub- and super-solution one gets on a point where $\sup ( u-\gamma v)$ is achieved , (see \cite{BD1} for    a complete proof using test functions) : 
                      $$ u^p -\mu u^{1+\alpha} r^{-2-\alpha} \leq | \nabla u |^\alpha F( D^2 u)\leq | \nabla ( \gamma v)|^\alpha F( D^2 \gamma v) \leq \gamma^{1+\alpha} v^p-\mu \gamma^{1+\alpha } v^{1+\alpha} r^{-2-\alpha}.$$
                       This contradicts the assumption on $\gamma$. 
                       
                       \end{proof}

    \begin{theo}\label{comp2}
    
     Suppose that $u$ and $v$ are respectively positive radial  sub-and super-solutions  and that there exist some constants $c_i$, $i = 1,2$ so that 
     $$c_1 \leq v^{p-1-\alpha} r^{2+\alpha}\leq c_2,  \ {\rm and } \ u^{p-1-\alpha} r^{2+\alpha} \leq c_2.$$
    Suppose  in addition that $v$ satisfies   $|v^\prime |\leq c r^{-(p+1) \over p-1-\alpha}$. 
            Then if $u\leq v$ on $r_o$, $u\leq v$ in all $B(0, r_o)\setminus \{0\}$.
      
      \end{theo}

    \begin{proof}

      Suppose that $\sup { u \over v} = \gamma^\prime >1$. Let $\gamma$ which is so that $\gamma >1$ and $\gamma^\prime$ first so that 
      $$\gamma < {1+\gamma^\prime \over 2}$$ 
      and $\gamma$ will be close to $\gamma^\prime$, chosen later. The supremum of $u-\gamma v$ is $>0$. 
            Suppose  that the supremum is achieved  on $\bar r \neq 0$.  In that case we easily reach the contradiction $\gamma^\prime \leq 1$ by arguing as in the proof of Theorem \ref{comp1}. 
      
      There remains to do the case where there exists some sequence coming to zero, so that $(u-\gamma v )(r_n) $ converges to the supremum. In particular $(u-\gamma v)(r_n)>0$. We claim that there exits 
      $n_1$ so that $u(r) > \gamma v(r) $ for any $r< r_{n_1}$. If not  for all $r_1$, there exist $\sigma < s < r_1$ so that $(u-\gamma v)(\sigma) = (u-\gamma v)(s) = 0$ and $u-\gamma v>0$ inside.  On the other hand we have 
      \begin{eqnarray*}
       {\cal M}^+ ( | \nabla u |^\alpha D^2 u)-| \nabla ( \gamma v)|^\alpha D^2 ( \gamma v)) &\geq & | \nabla u |^\alpha { \cal M}^+( D^2 u)-| \gamma \nabla v|^\alpha { \cal M}^+ ( D^2 ( \gamma v))\\
       &\geq & u^p-\gamma^{1+\alpha } v^p -\mu {u^{1+\alpha}-\gamma^{1+\alpha} v^{1+\alpha} \over r^{2+\alpha}}\\
       &\geq & 
        ( \gamma^p-\gamma^{1+\alpha}) v^p -\mu { ( \gamma^\prime)^{1+\alpha}-\gamma^{1+\alpha} \over r^{2+\alpha}} v^{1+\alpha} \\
        &\geq &( ( \gamma^p-\gamma^{1+\alpha}) c_1^p-\mu (( \gamma^\prime )^{1+\alpha}-\gamma^{1+\alpha}) c_2^{1+\alpha}) r^{-(2+\alpha)p\over p-(1+\alpha)}
        \end{eqnarray*}
        
        We impose to $\gamma$ to be so that 
        $$\mu \left(( {\gamma^\prime\over \gamma } )^{1+\alpha}-1\right) c_2^{1+\alpha}< {\left(\left({ \gamma^\prime + 1 \over 2} \right)^{p-1-\alpha}-1 \right)c_1^p\over 2}$$
        
      We have obtained  that  such $\sigma$ and $s$ do not exist since if they are, there exists some $t\in ]\sigma, s[$ so that 
        $(u^\prime -\gamma v^\prime) (t) = 0$ and 
        $$ (| u^\prime |^\alpha u^\prime -\gamma^{1+\alpha} |v^\prime |^\alpha v^\prime )^\prime (t) \leq 0.$$
         This contradicts the previous computations, in particular $u-\gamma v >0$ for $r< r_{n_1}$. We also get that there exists $c_3>0$ so that for $r$ small enough

        $$  {\cal M}^+ ( | \nabla u |^\alpha D^2 u-| \gamma \nabla v |^\alpha D^2( \gamma v)) \geq c _3r^{ -p(2+\alpha) \over  p-1-\alpha}= c_3 r^{-\beta p}$$
        
         with 
         $c_3 = {( \gamma^p-\gamma^{1+\alpha})c_1^p\over 2} $.
       Arguing as in the proof of   Proposition \ref{M+sign},  $u^\prime -\gamma v^\prime$ has a constant sign in a neighborhood of zero. Since 
         $$\lim_{n\rightarrow+\infty} u(r_n)-\gamma v(r_n) = \lim v(r_n) \left({u(r_n)\over v(r_n)}-\gamma\right) = +\infty,$$this sign is $\leq 0$.   One has 
         \begin{equation}\label{uprimev} c  r^{ {-p(2+\alpha) \over  p-1-\alpha}+ ( N^+-1)(1+\alpha)} \leq ((|u^\prime |^\alpha u^\prime -\gamma^{1+\alpha}|v^\prime |^\alpha v^\prime)  r^{ (N^+-1)(1+\alpha)})^\prime .
         \end{equation}
                             
        Integrating between $r$ and $2r$  the inequality (\ref{uprimev}), and using  
                 $ {2^k-1\over k} >0$
for $k = (N^+-1)(1+\alpha) +1-p\beta $ when it is $\neq 0$,          and also using $|u^\prime |^\alpha u^\prime (2r)-\gamma^{1+\alpha} |v^\prime |^\alpha v^\prime (2r) \leq 0$, we obtain that  : 
          
          If $k\neq 0$ 
          
           $$(\gamma ^{1+\alpha} |v^\prime |^\alpha v^\prime -  |u^\prime |^\alpha u^\prime) (r)\geq c r^{-(p+1)(1+\alpha) \over p-1-\alpha}= c r^{-(\beta+1)p}.$$
           
            While if $k = 0$, ( ie $N^+-2= {2+\alpha \over p-1-\alpha}$),   we get 
          $$ ( \gamma ^{1+\alpha} |v^\prime |^\alpha v^\prime -  |u^\prime |^\alpha u^\prime )(r)\geq c \log 2 r^{-(N^+-1)(1+\alpha)}= c \log 2  r^{-(p+1)(1+\alpha) \over p-1-\alpha}.$$
                        
            When $\alpha \leq 0$ we use   , depending on the signs of $u^\prime$ and $v^\prime$ the inequalities
             when $x, y, x-y>0$
             $$x^{1+\alpha} \leq y^{1+\alpha} + ( x-y)^{1+\alpha}$$ 
             
              and when $x>0, y<0$
              $$ x^{1+\alpha} + |y|^{1+\alpha} \leq 2^{-\alpha} ( x+ |y|)^{1+\alpha}$$
lead to                $$\gamma v^\prime -u^\prime \geq 2^\alpha c  r^{-(p+1) \over p-1-\alpha}.$$   
               
               In the case $\alpha >0$      
             we use the  additional assumption on $v^\prime $  and the fundamental calculus Theorem to get 
       $$ (\gamma^{1+\alpha} |v^\prime |^\alpha v^\prime - |u^\prime |^\alpha u^\prime )(r)\leq (1+\alpha)\sup_{]0, r[}  |\gamma v^\prime |^\alpha ( \gamma v^\prime -u^\prime(r) )\leq c r^{-(p+1)\alpha\over p-1-\alpha} ( \gamma v^\prime -u^\prime)(r).$$
       We derive   from this  that $$( \gamma v^\prime -u^\prime)(r) \geq c r^{-(p+1 )(-\alpha) \over p-(1+\alpha)}( \gamma^{1+\alpha} |v^\prime |^\alpha v^\prime - |u^\prime |^\alpha u^\prime )(r)\geq c r^{-(p+1 )(-\alpha) \over p-(1+\alpha)}  r^{-(p+1)(1+\alpha) \over p-1-\alpha}
       =cr^{-(p+1 )\over p-(1+\alpha)} ,$$
     
       in both cases  $\alpha >0$ and $\leq 0$. Integrating between $r$ and $2r$  sufficiently smalls once more,  one gets  that for some constant $c$ depending only on $\alpha$
      $$ (u-\gamma v) (r)\geq cc_3 r^{- \beta }$$

      always with some constant independant of $\gamma^\prime -\gamma$.      On the other hand one has by the definition of $\gamma$ and $\gamma^\prime $
      $$c  {( \gamma^p-\gamma^{1+\alpha})c_1^p\over 2} r^\beta \leq  (u-\gamma v )(r)\leq ( \gamma^\prime -\gamma) v(r)\leq c_2 ( \gamma^\prime -\gamma) r^{-\beta } .$$
      
       We have obtained a contradiction, since $\gamma$ and $\gamma^\prime$ can be taken arbitrarily close.
       \end{proof}
      
      \section{Proof of Theorem \ref{exi} : Existence of  solutions with the desired asymptotic behaviour}

       We begin to observe that if  $z = c r^{-\gamma}$,  with $\gamma >0$,  $z$ satisfies 
    $$ |z^\prime |^\alpha { \cal M}^+ ( D^2 z) + \mu { z^{1+\alpha} \over r^{2+\alpha}} = c^{1+\alpha} \left( \gamma ^{2+\alpha}-( N^+-2) \gamma^{1+\alpha} + { \mu \over \Lambda}\right) r^{-(\gamma+1)\alpha -\gamma-2}$$
    
     In particular,  if $\gamma \leq  \tau^-$ or $\gamma \geq \tau^+$ 
     $$|z^\prime |^\alpha { \cal M}^+ ( D^2 z) + \mu { z^{1+\alpha} \over r^{2+\alpha}} \geq 0, $$
     while if $\gamma \in ]\tau^-, \tau^+[$ 
    $$ |z^\prime |^\alpha { \cal M}^+ ( D^2 z) + \mu { z^{1+\alpha} \over r^{2+\alpha}} < 0.$$
    
     This implies that if $\gamma = \beta$, if $(\beta -\tau^-)(\beta-\tau^+)>0$ 
      and $c \leq K$, ( respectively $c> K$) $z$ is a sub-solution,  (respectively a super-solution).  For $c= K$,  $z$ is an exact solution. This proves the existence of solution behaving as $K r^{-\beta}$ in the cases $p< p^\star$ or $p> p^{\star\star}$.

  We are now  interested in the existence of solutions which satisfy  $\liminf u(r) r^\beta =0$.

     \begin{prop}\label{tau+}
           Suppose that $p < 1+\alpha + { 2+\alpha \over \tau^+}$,    Then for $A>0$ , and $a>0$ there exists   $\delta>0 $ small enough,   and $r_o(a, A,  \alpha, \tau^+, N^+, \delta)>0$ so that in $r< r_o$
            $\underline{w} = A r^{-\tau^+}- a r^{-\tau^++\delta}$ is a sub-solution of the equation. 
            
              Suppose that $p < 1+\alpha + { 2+\alpha \over \tau^-}$,   then for $A>0$, and for $a>0$,  there exists  $\delta>0$ small enough and $r_o:= r_o( a, A, \delta, \tau^-, N^+)  >0$ so that in $r< r_o$, 
            $\underline{w} = A r^{-\tau^-}+a r^{-\tau^-+\delta}$ is a sub-solution of the equation. 

            \end{prop}

  \begin{proof}

    We will use the inequalities :      For $\kappa = \alpha$ or $1+\alpha$, there exists some positive constant $C(\kappa)$ so that 
for any $u\in ]0, {1\over 2}[$,     $$1-\kappa u -C( \kappa ) u^2 \leq (1-u)^\alpha \leq 1-\kappa u + C( \kappa ) u^2$$

              We denote for simplicity $\tau_1 = \tau^+-\delta$. We choose $\delta $ small enough in order that  $\delta < ( 2+\alpha) ( \tau^+-\tau)$ and 
             $p \tau^+<       \tau^+(1+ \alpha )-\delta +\alpha+2    $.     We suppose first that $r_o^\delta < {A\tau^+\over 2(\tau_1+1) ^2 a }$, and $r<  r_o$.   Then  $\underline{w}$ is positive, decreasing and convex.              Having in mind that $\mu = -\Lambda(\tau^+)^{1+\alpha}  ( \tau^+-( N^+-2))$   
            and 
            $$ \tau_1 ( \tau_1-(N^+-2)) = -{\mu \over \Lambda( \tau^+)^\alpha} -{\delta} ( 2\tau^+-(N^+-2))+ { \delta^2 }$$
          One has 
             $$\underline{w} ^\prime = -\tau^+ A r^{-\tau^+-1} ( 1-{a \tau_1\over A \tau^+} r^\delta)$$
             
             $$  |\underline{w} ^\prime |^\alpha \leq     (\tau^+ A )^\alpha r^{-(\tau^++1)\alpha} (  1-{ \alpha a \over A} ( 1-{\delta \over \tau^+}) r^\delta +C( \alpha) \left( {a \tau_1\over A \tau^+} r^\delta)\right)^2) $$
             
             $$\underline{w}^{\prime \prime}+ \underline{w}^\prime {N^+-1\over r} =  A r^{-\tau^+-2}\left( -{\mu \over (\tau^+)^\alpha \Lambda} + {a \over A} \left( {\mu \over (\tau^+)^\alpha\Lambda} +\delta ({2\tau^+-(N^+-2)} )- {\delta^2 }) \right)r^\delta\right)$$

                 $$\underline{w}^{1+\alpha} \geq (Ar^{-\tau^+})^{1+\alpha}(1-{a\over A} (1+\alpha) r^\delta-C( 1+\alpha) \left({a\over A} r^\delta\right)^2)$$
                  From this one obtains that there exist some constants $c_1( a, A,  \alpha, \tau^+, N^+, \delta) $  and $c_2(a, A,  \alpha, \tau^+, N^+, \delta)  >0$,  so that 
                  \begin{eqnarray*}
                   |\underline{w} ^\prime |^\alpha ( \underline{w}^{\prime \prime}+ \underline{w}^\prime {N^+-1\over r} )&+& \mu {\underline{w} ^{1+\alpha} \over \Lambda r^{2+\alpha}}\\
                  &\geq & \delta A^\alpha a  (\tau^+)^{\alpha} ( \tau^+ ( 2+\alpha) -(1+\alpha) ( N^+-2)-\delta )  r^{-\tau^+ (1+\alpha)+\delta-\alpha-2}\\
                  &-& c_1( a, A,  \alpha, \tau^+, N^+, \delta) r^{-\tau^+ (1+\alpha)-\alpha-2 +2\delta}\\
                   &\geq : &c_2( a, A,  \alpha, \tau^+, N^+, \delta) r^{-\tau^+ (1+\alpha) +\delta -\alpha-2}
                   \end{eqnarray*} 
                    as soon as $r< r_1$, where $r_1$  is chosen small enough  depending on $a, A, \alpha, \tau^+, \delta,  N^+$. 
                    On another hand ,                    $$(\underline{w} )^p \geq  A^p  r^{-p\tau^+}.$$
                Choosing  $r_2$  so that  $r_2^{\tau^+(p-(1+\alpha))-\alpha -2+\delta} >{A^p\over c_2}$ , $r_2< r_1$, and $r< r_2$ 
                   we get the desired result.

        For the second assertion,      
    it is sufficient to replace in the calculation above $\tau^+$ by $\tau^-$ and $a$ by $-a$.  One gets   that there exists some constant $c_1(a, A,  \alpha, \tau^-, N^+, \delta) $, $c_2(a, A,  \alpha, \tau^-, N^+, \delta) >0$ , $c_3 ( a, A, p, \delta) 
>0$ so that 
              
           \begin{eqnarray*}
           |\nabla \underline{w}  |^\alpha { \cal M}^+ ( D^2 \underline{w}  )+ \mu { \underline{w} ^{1+\alpha }\over \Lambda r^{\alpha +2}}& =& -\delta a A^\alpha  ( \tau^-)^\alpha ( \tau^- ( 2+\alpha) -(1+\alpha) ( N^+-2)-\delta )  r^{-\tau^-(1+ \alpha) +\delta -\alpha-2}\\
           &-& c_1(a, A, \tau^-, \alpha, N^+, \delta)(r^{-\tau^- (1+\alpha) +2\delta -\alpha-2} )\\
       &          \geq &  c_2( a, A, \tau^-, \alpha, N^+, \delta )  r^{-\tau^+ (1+\alpha) +\delta -\alpha-2}\\
       &\geq &  c_3( a, A, p, \delta ) r^{-p\tau^-}
       \\
       &\geq &\underline{w} ^p
           \end{eqnarray*}
using  $ \tau^- ( 2+\alpha) -(1+\alpha) ( N^+-2) <0$,   as soon as $\delta$ is small enough in order that $p \tau^-<       \tau^-(1+ \alpha)-\delta +\alpha+2 $, and choosing $r_o$  conveniently. 
     \end{proof}

           \begin{proof} of Theorem \ref{exi}
       
       We provide a solution by  using a  decreasing sequence $(r_n)$ coming to zero, and  constructing a sequence $v_n$ on $[r_n, r_o]$ in the following way : 
             Let $r_o$ be small enough  and  $\underline{v}$ and $\overline {v}$ respectively  ordered sub- and super-solutions in $]0, r_o[$. 
        The first element of the sequence is  $\underline{v}$. We define $v_1$ on $[r_1,r_o]$ as the radial  solution of the equation 
        $$\left\{ \begin{array}{lc}
        -| \nabla v_1|^\alpha {\cal M}^+(D^2 v_1) + |v_1|^{p-1} v_1 = \mu { \underline{v}^{1+\alpha} \over r^{2+\alpha}}& \ {\rm on } \ [r_1,r_o[\\
        v_1(r_o) = \underline{v} (r_o)& v_1(r_1) = \underline{v} (r_1)\end{array}\right.$$
        
         The existence of $v_1$ is a consequence of the existence's results in \cite{BD1}, \cite{BD2}. Furthermore, by the comparison principle in \cite{BD1}, $v_1 \geq \underline{v} > 0$ on $[r_1, r_o]$. Once $v_n$ is defined on $[r_n, r_o]$ we define $v_{n+1}$ on $[r_{n+1}, r_o]$ by
           $$\left\{ \begin{array}{lc}
        -| \nabla v_{n+1}|^\alpha {\cal M}^+(D^2 v_{n+1}) + |v_{n+1}|^{p-1} v_{n+1} &=\left\{\begin{array}{cc}\mu {v_n^{1+\alpha} \over r^{2+\alpha} } & {\rm on } \ [r_n, r_o[\\
         \mu { \underline{v}^{1+\alpha} \over r^{2+\alpha}}& \ {\rm on } \ [r_{n+1},r_n[
         \end{array}\right.\\
        v_{n+1}(r_o) = \underline{v} (r_o)  , \ v_{n+1}(r_{n+1}) =  v_n(r_{n+1}) &
        \end{array}\right.$$
$v_{n+1}$ is well defined and positive. By the classical  comparison  principle in \cite{BD1}, $(v_n)_n$ is increasing, and 
$\underline{v} \leq v_n \leq \overline{v}$. on $[r_n, r_o]$.  Standard stability and regularity results in  \cite{BD1} then imply that the monotone locally uniformly bounded sequence $\{v_n\}$ converges in $C^2_{{\rm loc}}((0,r_o])$ to a radial solution $v$ of equation \eqref{eqmain}. Furthermore using Proposition \ref{tau+} and choosing $A< B$,  and  $a>0$, for  $r_o$ convenient, 
$\overline{w} = B r^{-\tau^+}$ and 
$\underline{w} = A r^{-\tau^+}- a r^{-\tau^++\delta}$ are respectively super- and sub-solutions  of the equation, on $]0, r_o[$, so that $\underline{w} < \overline{w}$.
 We get that we can construct a solution such that  for $r$ small enough and for some positive constants $c_i$, $i=1;2$ : 
 $$c_1  \leq v r^{\tau^+} \leq c_2$$
It is  left to prove that $v(r)$ is asymptotic to $r^{-\tau^+}$ as $r\to 0$. 
By the bounds above,  there exists   a decreasing sequence $\{ r_n\}$ converging  to zero  such that  $v(r_n) r_n^{\tau^+} \rightarrow c_1$ for some $c_1>0$. 
For $\epsilon >0$   arbitrarily small, and $\delta$ as in Proposition \ref{tau+},   let $\nu \geq 1$ such that
  $ c_1-\epsilon \leq  v(r_n) r_n^{\tau^+} \leq c_1+\epsilon$ for $n\geq \nu$ and, moreover,  $r_\nu<r_o(\epsilon, c_1-\epsilon, \tau^+, N^+, \delta)$ as in Proposition \ref{tau+} 
  
Then, for all $n\geq \nu$, we can consider on the interval $[r_{n+1}, r_n]$ the supersolution
     $\overline{z}  = (c_1+\epsilon) r^{-\tau^+}$ and the subsolution 
     $\underline{z}= (c_1-\epsilon)  r^{-\tau^+}-\epsilon   r^{-\tau^++\delta }$. 
 Since     
 $$ \underline{z}(r_n) \leq (c_1-\epsilon ) r_n^{-\tau^+} \leq v(r_n)\leq (c_1+\epsilon) r_n^{-\tau^+}= \overline{z}(r_n),$$ and 
$$ \underline{z}(r_{n+1} )\leq (c_1-\epsilon ) r_{n+1}^{-\tau^+} \leq v(r_{n+1})\leq (c_1+\epsilon) r_{n+1}^{-\tau^+}= \overline{z}(r_{n+1})\, ,$$
       the comparison principle given by Proposition \ref{comp1} yields 
       $\underline{z}\leq v \leq  \overline{z} $ in $[r_{n+1}, r_n]$ for all $n\geq \nu$, which implies  that finally 
       $$c_1-2\epsilon \leq v(r) r^{\tau^+} \leq c_1+\epsilon$$
        
\medskip                                          

       If $p< p^{\star\star}$, in order  to prove the existence of a solution equivalent to some constant $r^{-\tau^-}$ , 
      we use the same construction as in the first part concerning $\tau^+$.  We obtain a solution which satisfies 
       $c_1 \leq v r^{\tau^-} \leq c_2$. To prove that $v$ is  asymptotic to $r^{-\tau^-}$,
      by the bounds above,  there exists   a decreasing sequence $\{ r_n\}$ converging  to zero  such that  $v(r_n) r_n^{\tau^+} \rightarrow c_1$ for some $c_1>0$. 
For $\epsilon >0$ arbitrarily small, let $\nu \geq 1$ such that
  $ c_1-\epsilon \leq  v(r_n) r_n^{\tau^-} \leq c_1+\epsilon$ for $n\geq \nu$ and, moreover,  $r_\nu< r_o( \epsilon, c_1-2\epsilon,\delta, \tau^-, N^+)$ as in Proposition \ref{tau+}.  
Then, for all $n\geq \nu$, we can consider on the interval $[r_{n+1}, r_n]$ the supersolution 
     $\overline{z}  = (c_1+\epsilon) r^{-\tau^-}$ and the subsolution 
     $\underline{z}= (c_1-2\epsilon)  r^{-\tau^-}+ \epsilon   r^{-\tau^-+\delta }$.     One has 
      $\underline{z}\leq \overline{z} $. 
 Since     
 $$ \underline{z}(r_n) \leq (c_1-\epsilon ) r_n^{-\tau^-} \leq v(r_n)\leq (c_1+\epsilon) r_n^{-\tau^+}= \overline{z}(r_n)$$
 as well as 
     $$ \underline{z}(r_{n+1} )\leq (c_1-\epsilon ) r_{n+1}^{-\tau^-} \leq v(r_{n+1})\leq (c_1+\epsilon) r_{n+1}^{-\tau^-}= \overline{z}(r_{n+1})\, ,$$
       the comparison principle given by Theorem  \ref{comp1} yields 
       $\underline{z}\leq v \leq  \overline{z} $ in $[r_{n+1}, r_n]$ for all $n\geq \nu$, which implies  that finally 
       $$c_1-2\epsilon \leq v(r) r^{\tau^-} \leq c_1+\epsilon$$
       
        For the case $ p = p^{\star\star}$ it is sufficient to remark that 
        for $a < \tilde K$, and  for  any $c>0$ 
        $$ \underline{w} = a r^{-\tau^-} ( -\log r+ c)^{-\tau^-\over 2+\alpha}$$
        is a sub-solution for $r$ small enough, 
         while if $a > \tilde K$ it is a super-solution. 
         We do not give the details of the proof  : We just remark that  for $a< \tilde K$
         \begin{eqnarray*}
         | \underline{w}^\prime|^\alpha (   \underline{w}^{\prime\prime} + {N^+-1\over r}  \underline{w}^\prime) + \mu  \underline{w}^{1+\alpha} r^{-2-\alpha} 
         &=&(a\tau^-)^\alpha   ( a\tau^-(\tau-\tau^-)) a^{-p} \underline{w} ^p+  O( r^{-\tau^-p} (c -\log r)^{-\tau^-p-1})\\
         &=&  \Lambda ^{-1}a^{1+\alpha-p} \tilde K^{p-1-\alpha}  \underline{w} ^p+ o( \underline{w} ^p) \\
         &\geq& \Lambda ^{-1}   \underline{w} ^p.
         \end{eqnarray*}
         as soon as $r$ is small enough. 
         An analogous calculation permits to prove that for $a>\tilde K$, $\overline{w} =  a r^{-\tau^-} ( -\log r+ c)^{-\tau^-\over 2+\alpha}$ is a super-solution,  as soon as $r$ is small enough. 
          The  reasoning for the end of the proof is as in the case $p< p^{\star\star}$.

         \end{proof}

  \bigskip

         \section{Proof of Theorem \ref{behav}}
         
       We begin with  some preliminary observations.

First of all, we note that if $u$ is a nonnegative solution of equation \eqref{eqmain}, by the equation in radial coordinates, one gets that on a point $r_o\neq 0$,  if $u(r_o) = 0$, then $u^\prime (r_o) = 0$ since $r_o$ is a minimum and then by Cauchy Lipschitz Theorem, $u$  vanishes identically or it is strictly positive in $B_1\setminus \{0\}$. Thus, in the sequel, by a solution  we mean a positive radial solution.

\begin{prop}\label{M+sign}

Suppose that $| \nabla u |^\alpha {\cal M}^+( D^2 u)$  has a sign for $r$ small enough, and $u$ is unbounded. Then $u^\prime \leq 0$
\end{prop}

 \begin{proof} 
  If $| \nabla u |^\alpha { \cal M}^+( D^2 u) \geq 0$, since $u$ is radial, this can be written as 
  $$ \Lambda\left( {(|u^\prime |^\alpha u^\prime ) ^\prime\over 1+\alpha}\right)^+-\lambda \left( {(|u^\prime |^\alpha u^\prime ) ^\prime\over 1+\alpha}\right)^- + \Lambda {N-1\over r}   (|u^\prime |^\alpha u^\prime )^+-\lambda  {N-1\over r} (|u^\prime |^\alpha u^\prime )^-\geq 0$$
  Suppose that $u^\prime$ changes sign indefinitely, then in particular there exists $(r_n)_n$ coming to zero  so that 
 $ u^\prime (r_n) = 0$ and $u^\prime <0$ in $]r_{n+1}, r_n[$. By the mean value's theorem there exists $s_n \in ]r_{n+1}, r_n[$ so that 
 $u^{\prime \prime} (s_n) = 0$,  a contradiction since  by the equation 
 $ \Lambda |u^\prime |^\alpha u^{\prime \prime}(s_n) > 0$.  
  The reasoning is analogous if $| \nabla u |^\alpha { \cal M}^+(D^2u) \leq 0$. 
   Then $u^\prime$ has constant sign near zero, but if $u^\prime \geq 0$, this would  contradict $u$ unbounded. 
   \end{proof}

We  prove that any non trivial positive radial solution satisfies  :

\begin{prop}There exists some constant $C$ so that 
$$u(r) \leq C r^{-2-\alpha \over p-1-\alpha}= C r^{-\beta}$$

\end{prop}

\begin{proof}

 We preliminarily observe that if $u(r)r^{\frac{2+\alpha}{p-1-\alpha}}$ is bounded along a decreasing sequence $\{r_n\}$ converging to zero, then it is bounded in a neighborhood of zero. Indeed, let us assume that there exists $C>0$ such that 
$$
u(r_n)r_n^{-\beta}\leq C\, .
$$
By using the computation made at the beginning of the section, we know that, for $C>0$ sufficiently large, the function 
 $C r^{-2-\alpha\over p-1-\alpha}$ is a super-solution of  equation \eqref{eqmain}. Then, by applying Theorem 
 \ref{comp1}   in every  annulus  $B_{r_{n}}\setminus\overline{B}_{r_{n+1}}$, we obtain  
$$u(r)\leq C r^{-\beta} \quad \hbox{ for } 0<r\leq r_1\, .
$$ 
Hence, if we assume  that $u(r)r^{\beta}$ is unbounded along a decreasing sequence $\{r_n\}$ converging to zero then we deduce
$$
\lim_{r\to 0} u(r) r^\beta = +\infty\, .
$$    
This implies that, for $r$ small enough, one has
$$
| \nabla u |^\alpha \mathcal{M}^+(D^2u)= u^p -\mu \frac{u^{1+\alpha}}{r^{2+\alpha}} \geq { u^p \over 2}>0\, .$$ 

  Then,  by Proposition \ref{M+sign}, 
 since $u$ is unbounded, $u^\prime \leq 0$.  Then $u^{\prime \prime} \geq 0$ and 
 by the equation there exists some positive constant $c$ so that  $|u^\prime |^\alpha u^{ \prime \prime} \geq cu^p.  $
   Multiplying by $u^\prime \leq 0$ one gets  for another positive constant $c$
   $$ ( |u^\prime |^{\alpha+2}  -{ c }u^{ p+1})^\prime \leq 0, $$
    integrating one gets  for $r< r_o$ 
   $$  |u^\prime |^{\alpha+2}  (r_o) - {c} u^{p+1} (r_o) \leq  |u^\prime |^{\alpha+2} (r) - {c} u^{p+1}( r) .$$
    Using $u^{p+1} \rightarrow +\infty$ , there exists some constant $c$ which can vary from one line to another so that for $r_o$ small enough and $r< r_o$, 
    $$ {|u^\prime |^{\alpha +2}\over  u^{p+1} } \geq c, $$ which means 
    $$ { u^\prime \over u^{ p+1 \over 2+\alpha} } \leq -c.$$ 
     Integrating one gets  using $\lim _{r\rightarrow 0} u^{1-{ p+1 \over 2+\alpha}}= 0$ that 
     $$u^{1-{ p+1 \over 2+\alpha}}\geq c r,$$
     but this contradicts the assumption. 
     
     \end{proof} 
     
     The following estimate on $u^\prime$ will be needed for some of the results. It  also  permits to obtain the results in most of the cases by an alternative proof,  with the aid of the Theorem of Poincar\'e Bendixson.

     \begin{prop}\label{uprime}
 Suppose that $\beta \neq N^+-2$.  Suppose that $u$ is a positive radial solution, so that $u^\prime \leq 0$.  There exists some constant $C$ so that 
  $$ |u^\prime | \leq C r^{-\beta-1}.$$
   \end{prop}

    \begin{proof}
   
    Let us observe that  whatever is the sign of $u^{\prime \prime}$
    $$|u^\prime |^\alpha ( \Lambda u^{\prime \prime} + \lambda {(N-1) u^\prime \over r}) \leq u^p-\mu u ^{1+\alpha} r^{-2-\alpha} $$
    
    Then 
     \begin{equation}\label{int}(|u^\prime |^\alpha u^\prime  r^{(N^+-1) (1+\alpha)})^\prime \leq c r^{-(\beta+1)(1+\alpha)-1+ (N^+-1)(1+\alpha)}.
     \end{equation}
    On  the other hand,  if $u^{\prime \prime} \leq 0$ by the equation, 
    
     $$|u^\prime |^\alpha( u^{\prime\prime} + { \lambda \over \Lambda} {(N-1) u^\prime \over r}) \geq  |u^\prime |^\alpha( u^{\prime\prime} +  {(N-1) u^\prime \over r})= { 1\over \lambda} (u^p-\mu u ^{1+\alpha} r^{-2-\alpha} )$$
        and 
     if      $u^{\prime \prime} \geq 0$ 
    $$| u^\prime |^\alpha (  u^{\prime \prime} +  {(N^+-1) u^\prime \over r}) = {1\over \Lambda} ( u^p-\mu u ^{1+\alpha} r^{-2-\alpha} )$$
      Then  for some positive constant $c$
$$   |\nabla u |^\alpha { \cal M}^+( D^2u) \geq   \inf ( { 1\over \lambda} (u^p-\mu u ^{1+\alpha} r^{-2-\alpha} ), { 1\over \Lambda} (u^p-\mu u ^{1+\alpha} r^{-2-\alpha} )) \geq -c r^{-\beta p}$$
and then
   \begin{equation}\label{int1}(|u^\prime |^\alpha u^\prime  r^{(N^+-1) (1+\alpha)})^\prime \geq -c r^{-(\beta+1)(1+\alpha)-1+ (N^+-1)(1+\alpha)}.
     \end{equation}

     Suppose first that $\beta > N^+-2$,  :    Integrating  \eqref{int} between $r$ and $r_o> r$ one gets 
      $$ (|u^\prime |^\alpha u^\prime) (r_o) r_o^{(N^+-1) (1+\alpha)}- (|u^\prime |^\alpha u^\prime)(r)   r^{(N^+-1) (1+\alpha)}\leq {-c\over \beta-N^++2}(r_o^{(N^+-2-\beta)(1+\alpha)}-r^{(N^+-2-\beta)(1+\alpha)} )$$
       and then  for some constants $c_1$,  and $c_2>0$
     $$|u^\prime |^\alpha u^\prime(r)  \geq c_1 r^{(N^+-1)(1+\alpha)} -c_2 r^{-(\beta+1 )(1+\alpha)}\geq -c_3 r^{-(\beta+1 )(1+\alpha)}$$

     If $\beta < N^+-2$ , \eqref{int1} 
     implies that $|u^\prime |^\alpha u^\prime  r^{(N^+-1) (1+\alpha)}$ has a limit since 
      $$(|u^\prime |^\alpha u^\prime  r^{(N^+-1) (1+\alpha)}+ c r^{-(\beta+1)(1+\alpha)+ (N^+-1)(1+\alpha)})^\prime \geq 0$$
       and $r^{-(\beta+1)(1+\alpha)+ (N^+-1)(1+\alpha)}$ tends to zero by the assumption on $\beta$. Note that if $|u^\prime |^\alpha u^\prime  r^{(N^+-1) (1+\alpha)}$ does not tend to zero, one has 
       for some constant $c>0$
        $$|u^\prime |^\alpha u^\prime  r^{(N^+-1) (1+\alpha)}\leq -c$$
         which implies that $u^\prime \leq -c r^{1-N^+}$ and then $u \geq c r^{2-N^+}$,  a contradiction  with 
         $u \leq c r^{-\beta}$.
          Then 
          $|u^\prime |^\alpha u^\prime  r^{(N^+-1) (1+\alpha)}$ tends to zero, and integrating \eqref{int1} between $0$ and $r$ one gets 
          $$| u^\prime |^\alpha u^\prime  r^{(N^+-1) (1+\alpha)} \geq -c_3 r^{-(\beta+1 )(1+\alpha)} r^{(N^+-1) (1+\alpha)}$$
           which implies $u^\prime \geq -c r^{-\beta-1}$.

            In particular in each of the cases,  $x(\log r) = r^\beta u(r)$ satisfies  $(x,x^\prime)$ is bounded. 
            \end{proof}

    \begin{prop}\label{signuprime}
    Suppose that $u$ is a positive radial solution of the equation. 
     Suppose that $x = u(r) r^\beta$ tends to $l\neq 0$, $l\neq    \mu ^{1\over p-1-\alpha}$, or tends to $  \mu ^{1\over p-1-\alpha}$ monotonically. Then  $| \nabla u |^\alpha {\cal M}^+( D^2 u)$  has a sign for $r$ small enough, and then $u^\prime \leq 0$.
     
     \end{prop}
      \begin{proof} The hypothesis implies $u$ unbounded. If $l \neq   \mu ^{1\over p-1-\alpha}$, $| \nabla u |^\alpha { \cal M}^+(D^2u) $ has the sign of $l^{p-1-\alpha}-\mu$ as soon as $r$ is small enough .  If $x\rightarrow  \mu ^{1\over p-1-\alpha}$ and $x^\prime \leq  0$ ( respt $<0$),  the right hand side has still a sign and we conclude as before. 
      \end{proof}

     \begin{theo}\label{useconde}
     
      Suppose that $u$ is a  non trivial   solution which satisfies 
      $$\lim u(r) r^\beta=0$$
       then 
       
       $$u^\prime \leq 0, \ { \rm and }  \    ( |u^\prime |^\alpha u^\prime )^\prime   \geq 0.$$
       \end{theo}
       
        \begin{proof}
        
        From the equation one gets for $r$ small enough 
         $$| \nabla u |^\alpha { \cal M}^+( D^2 u) \leq 0.$$
                By  Proposition \ref{M+sign}, $u^\prime \leq 0$.          We prove that $(|u^{ \prime}|^\alpha u^ \prime)^\prime $ does not change sign. Suppose that we have a sequence $r_n$ coming to zero, so that $(|u^\prime |^\alpha u^\prime )^\prime (r_n) = 0$, and $(|u^\prime |^\alpha u^\prime )^\prime >0$ inside the intervall. Then one must have $(|u^\prime |^\alpha u^\prime )^{\prime\prime}(r_{n+1} )>0$ while $(|u^\prime |^\alpha u^\prime )^{\prime\prime}(r_n) <0$.   
           
         To get a contradiction ,  it is sufficient to prove that 
           $$r\mapsto -\lambda ( N-1)r^{-1} |u^\prime |^\alpha(r) u^\prime(r) + u^p(r) -\mu u ^{1+\alpha}(r) r^{-2-\alpha} $$ has a derivative $> 0$ on a point where $(|u^\prime |^\alpha u^\prime )^\prime $ is zero. 
           
            We consider 
            $$\varphi(r) =  -\lambda (N-1)r^{1+\alpha} |u^\prime |^\alpha u^\prime +u^p r^{2+\alpha} -\mu u^{1+\alpha}.$$
            Then  $$\varphi^\prime (r) = - \lambda (N-1)) r^{1+\alpha} ( |u^\prime |^\alpha u^\prime)^\prime + p u^{p-1}  u^\prime r^{2+\alpha}  -\mu (1+\alpha) u^\alpha u^\prime - \lambda ( N-1) (1+\alpha) r^\alpha  |u^\prime |^\alpha u^\prime+ (2+\alpha) r^{1+\alpha} u^p$$
              has the sign of $- u^\alpha u^\prime\left( \mu (1+\alpha)-pu^{p-1-\alpha} r^{2+\alpha} \right)  - \lambda ( N-1) (1+\alpha) r^\alpha |u^\prime |^\alpha u^\prime + (2+\alpha) r^{1+\alpha} u^p$  on a point where 
              $(|u^\prime|^\alpha u^ \prime)^\prime  =0$, hence $> 0$, a contradiction on $r_n$. Since $u$ is unbounded, and  $(|u^\prime |^\alpha u^\prime)$ is monotone and negative, it tends to $-\infty$. 
              
            If we now suppose by contradiction that $(|u^\prime |^\alpha u^\prime )^\prime <0$, this leads  to a contradiction with the fundamental calculus Theorem. Indeed,  since $u^\prime \rightarrow -\infty$, there exists $r_o$ so that  for $r< r_o$
            $|u^\prime|^{1+\alpha} (r_o) < {|u^\prime |^{1+\alpha}(r) \over 2}$, then  
            $$0< {|u^\prime |^\alpha u^\prime (r_o)-|u^\prime |^\alpha u^\prime (r) \over r_o-r} = (|u^\prime |^\alpha u^\prime )^\prime (r_1)$$
            for some $r_1< r_o$. 
            We have obtained that $ (|u^\prime |^\alpha u^\prime )^\prime(r)\geq 0$ for $r$ small enough. 
              
              \end{proof}

\bigskip
              
              In order to precise the differential equation satisfied by $x$ in the sequel,   let us devote a few lines to  precise the region in the phase space $(x, x^\prime)$ which corresponds to the region where $(|u^\prime |^\alpha u^\prime )^\prime \geq 0
$.
 Note that 
 $$(|u^\prime |^\alpha u^\prime )^{ \prime} \geq 0 $$
  if and only if
   \begin{equation}\label{xuseconde}
    \left\{\begin{array}{cc} x^{p-1-\alpha} \geq  \mu, &\ \\
  
{\rm Or} \   x^{p-1-\alpha} \leq \mu,  \ {\rm  and } \ x^\prime\leq \beta x -( {\mu x^{1+\alpha} -x^p\over \Lambda ( N^+-1)})^{1\over \alpha+1}&
  \end{array}\right.
  \end{equation}
  
   Note for further purposes that if $p< p^{\star\star}$, $\beta^{1+\alpha} > \mu \Lambda^{-1}$ and then a point  $\bar t$ where $x^\prime=0$ is in the region where $u^{\prime \prime } (e^ t
)\geq 0$. Indeed, $p<  p^{\star\star}$ implies $\beta^{1+\alpha}  > (\tau^-)^{1+\alpha}  $
 and 
 
$$ {\mu \over \Lambda (N^+-1)} = (\tau^-)^{\alpha+1} (N^+-2-\tau^-)
(N^+-1)^{-1}\leq  (\tau^-)^{\alpha+1} $$   (This permits to prove in particular , -even if we will not need this here-  that there is no periodic orbit for the equation, for $t< T$ large enough negative  ).              
                                    
                                   \bigskip
                                   
                                        \subsection { The case  where lim$u(r) r^\beta=0$.   }        
         The following theorem proves the last assertions in item 1
     . 
        
         \begin{theo}\label{lim=0}
          Suppose that $p< 1+ \alpha +{2+\alpha \over \tau^+}$
          Then any positive solution $u$  of \eqref{eqmain} so that $\lim u(r) r^{2+\alpha \over p-1-\alpha}  = 0$ satisfies 
          $$u(r) \sim c r^{ -\tau^\pm}.$$ 
          
          \end{theo}
          
          The proof of Theorem \ref{lim=0} is the consequence of several propositions : 
          
          \begin{prop}\label{propy}
  
    Suppose that $p\leq p^{\star \star}$.  Let $u$ be a positive solution  of \eqref{eqmain}. Let $x( \log r) = r^{\beta} u( r)$ and $y = \beta-{x^\prime \over x}$. 
    If $x\rightarrow 0$ at $-\infty$
    there exists $T$ so that for $t< T$, $y^\prime (t) \neq 0$, and then $y$ has a limit when $t$ goes to $-\infty$. 
    
     \end{prop}
     \begin{proof}

      If $x\rightarrow 0$, we know that $u^{\prime\prime} > 0$, for $r$ small enough,  then we are on a region where $u(r) $, $x(\log r)$ and also then $y(\log r)$ are ${ \cal C}^\infty$. We also have $y\geq 0$, since $u^\prime \leq 0$. The equation satisfied by $y$ is 
      \begin{equation}\label{eqy} (y^{1+\alpha}  )^\prime ( 1+\alpha)^{-1} = y^{2+\alpha}-(N^+-2) y^{1+\alpha} + { \mu -x^{p-1-\alpha} \over \Lambda}.
      \end{equation}
      Note that 
               \begin{equation}\label{defpsi}
               \psi(y) =  y^{2+\alpha}-(N^+-2) y^{1+\alpha} + \Lambda^{-1}\mu
               \end{equation}
                satisfies $\psi ( \tau^\pm) = 0$ and is $>0$ in $]0, \tau^-[ \cup ]\tau^+, +\infty[$ and is $<0$ in $]\tau^-, \tau^+[$. 
                                We prove that $(y^{1+\alpha})^\prime$ has constant sign for $t$ large enough negative. 
                 Indeed if not there exists a sequence $t_n$ coming to $-\infty$, so that $(y^{1+\alpha})^\prime(t_n) = (y^{1+\alpha})^\prime(t_{n+1})=0$ and $(y^{1+\alpha})^{\prime }>0$ in $]t_{n+1}, t_n[$ .  
                 Let $\epsilon >0$ so that $|\psi(y) | < \epsilon$ implies $\inf ( |y-\tau^+|, |y-\tau^-|)  < { \beta-\tau^+\over 2}$ and let $T$ be large enough in order that $t< T$ implies $x^{p-1-\alpha} (t) <\Lambda  \epsilon$. Then if $T_1<T$ is so that $(y^{1+\alpha})^\prime (T_1) = 0$, by the equation 
                 $\psi(y(T_1)) ={ x^{p-1-\alpha}\over \Lambda} < {\epsilon }$, 
                  and then 
                  $\inf( |y(T_1)-\tau^+|, |y(T_1)-\tau^-|)  < { \beta-\tau^+\over 2}$,  implying ${x^\prime \over x} (T_1)> { \beta-\tau^+\over 2} >0$.  Then derivating the equation in $y$ one gets 
                  $$({y^{1+\alpha}\over 1+\alpha})^{\prime \prime}  (T_1)=(  y^{2+\alpha} -(N^+-2) y^{1+\alpha})^\prime(T_1) -(p-1-\alpha) x^{p-1-\alpha} {x^\prime\over x\Lambda}(T_1)  <0$$
                  This yields a contradiction when  $T_1=t_{n+1}$. We have obtained that $y^\prime$ does not change sign.
                  
                  \end{proof}
            The two following propositions  precise the exact behavior of $u$   when $x$ goes to zero :       
         \begin{prop}\label{phi}
         Suppose that $u$ is a positive solution of \eqref{eqmain}. Suppose that   $\tau^+< \beta$ and 
         $${ u^\prime r\over u} = -\tau^+ + \varphi, $$
          and that  for $r$ small enough, 
          $\varphi, r\varphi^\prime$ have the same constant sign , and  both  tend to zero.  (In particular $\lim u r^\beta = 0$). There exist positive constants $c_1^+, c_2^+, c_3^+$  and $r_o$ so that for $r< r_o$

          $${\rm if } \ \varphi, \varphi^\prime >0, \ r\varphi^\prime -c_1^+ \varphi \geq c_2^+ u^{p-1-\alpha} r^{2+\alpha}, \ 
                    {\rm if } \ \varphi, \varphi^\prime <0, \ r\varphi^\prime -c_1^+ \varphi \leq c_3^+ u^{p-1-\alpha} r^{2+\alpha}.$$

As a consequence there exists $c_4^+, c_5^+$  positive so that for $r$ small enough 
$$ |\varphi| \leq c_4^+ r^{c_5^+}.$$
More precisely, for all $\epsilon$ we can choose $r_o$ small enough in order that  the inequality above holds for $r< r_o$ with  
 $|c_5^+ -( \inf ((2+\alpha) ( \tau^+-\tau), 2+\alpha-\tau^+(p-1-\alpha))|  = O( \epsilon)) $, in particular $c_5^+>0$)  .

             Finally in each of the cases  there exists some constant $c>0$ so that 
             $$u(r) \sim c r^{-\tau^+}.$$
             
             \end{prop}
                       \begin{rema}
             The fact that $\varphi$ and $\varphi^\prime$ have  a  sign is not necessary, but Proposition \ref{phi} will be used   later  only  in these situations.               \end{rema}

              \begin{prop}\label{phi1}
               Suppose that $u$ is a positive solution of \eqref{eqmain}, and  that  $\tau^- < \beta$, and 
               $$ {u^\prime r \over u} = -\tau^-+\varphi, $$
                with $\varphi \geq 0$, $r\varphi^\prime \geq 0$  for $r$ small enough, and both tend to zero. 
                There exist positive constants $c_1^->0,  c_2^->0$ so that 
          
                     $$r\varphi^\prime +c_1^- \varphi \leq c_2^- u^{p-1-\alpha} r^{2+\alpha}$$
               
           As a consequence,   for all $\epsilon$,  there exists  $c_3^-$,  and $r_o$ so that 
for $r< r_o$ 
           $\varphi (r) \leq {c_3^-}  r^{-\tau^-(p-1-\alpha)+2+\alpha-\epsilon(p-1-\alpha)} $.

                         Finally one has, for some positive constant $c$  
           $$ u(r) \sim c r^{-\tau^-}.$$
           \end{prop}
             
             \begin{proof}of Proposition \ref{phi}
              $$ {u^{\prime \prime }r^2\over u} = {(u^\prime )^2 r^2 \over u^2}-{u^\prime r \over u} + \varphi^\prime r= ( \tau^+-\varphi)^2+ ( \tau^+-\varphi) + \varphi^\prime r$$
               and then

              \begin{eqnarray*}
                (\tau^+)^\alpha ( 1-\alpha { \varphi \over \tau^+}&+& O( \varphi ^2) )\left( r\varphi^\prime + ( N^+-2-2\tau^+)\varphi  +  \varphi^2\right.\\
               &+&\left. \tau^+( \tau^+-(N^+-2))\right) +\Lambda^{-1}  \mu  \\
            &=&  |{ u^\prime r \over u}|^\alpha((N^+-2-2\tau^+) \varphi + r\varphi^\prime + \tau^+ ( \tau^+-(N^+-2))+\varphi^2) + \Lambda^{-1}\mu \\
             & = &|u^\prime r u^{-1}|^\alpha({ u^{\prime\prime}r^2\over u} + ( N^+-1 ){ u^\prime r\over u} ) + \Lambda^{-1}\mu  
                 = \Lambda^{-1}u^{p-1-\alpha} r^{2+\alpha}
               \end{eqnarray*}
               
               If $\varphi \geq 0$, using $\varphi, r\varphi^\prime$ tend to zero and $\tau^+\leq N^+-2$,  
               for all $\epsilon >0$ , $\epsilon <\inf (  {1\over 4}, {(2+\alpha) ( \tau^+-\tau ) \over 4 ( N^+-\tau^+)})$ there exists $r_o$ so that for $r< r_o$, 
                \begin{equation}\label{epsilonphi}
                {1+|\alpha|\over \tau^+} ( | N^+-2-2\tau^+| + 1)  \varphi(r)< \epsilon
                \end{equation}
                 
                $$(1-(\alpha+\epsilon) { \varphi \over \tau^+})( r\varphi^\prime +\varphi (N^+-2-2\tau^+ +\epsilon)+\tau^+ ( \tau^+-(N^+-2)) )+ {\mu \over \Lambda (\tau^+)^\alpha}  \geq  { u^{p-1-\alpha} r^{2+\alpha }\over( \tau^+)^\alpha \Lambda}$$
               
                 Using the value of $\mu$ , the hypothesis on $\epsilon$  and $c_1 = -(1+\alpha) ( N^+-2)+(2+\alpha) \tau^+= (2+\alpha) ( \tau^+-\tau)>0$, there exists  $r_1< r_o$  so that for  $r< r_1$ 
              
             $$  r\varphi^\prime(1+\epsilon) -(c_1+2\epsilon ( \tau^+-N^+))  \varphi  \geq  {1\over \Lambda  (\tau^+)^\alpha }  u^{p-1-\alpha} r^{2+\alpha}
           . $$
We get immediately  with $c_1^+ := {(c_1+2\epsilon ( \tau^+-N^+))\over 1+\epsilon} $ 
$$(r^{-c_1^+} \varphi)^\prime \geq 0,$$
                which implies by integrating between $r$ and $r_o$
                
                $$ \varphi(r) \leq \varphi(r_o) r_o^{-c_1^+} r^{c_1^+}.$$
                                   
              If $\varphi, \varphi^\prime \leq 0$, for all $\epsilon >0$ there exists $r_o$ so that for $r< r_o$
                $$ (1-(\alpha+\epsilon) { \varphi \over \tau^+})( r\varphi^\prime  +\varphi (N^+-2-2\tau^+ +\epsilon)+\tau^+ ( \tau^+-(N^+-2)) )+ {\mu \over \Lambda (\tau^+)^\alpha} 
                 \leq { u^{p-1-\alpha} r^{2+\alpha }\over( \tau^+)^\alpha \Lambda}  $$
                            Using the value of $\mu$, choosing $\epsilon$ and $r_o$ as in  \eqref{epsilonphi},  and $c_1 = -(1+\alpha) ( N^+-2)+(2+\alpha) \tau^+= (2+\alpha) ( \tau^+-\tau)>0$ , there exists  $r_1< r_o$   so that for $r< r_1$ 
               \begin{equation}\label{varphi} (\tau^+)^\alpha ( r\varphi^\prime(1+\epsilon) -(c_1+2\epsilon ( \tau^+-N^+))  \varphi  )\leq  {1\over \Lambda}  u^{p-1-\alpha} r^{2+\alpha}. \end{equation}

                           On the other hand, by the definition of $\varphi$,
   fixing $\epsilon < \tau^+$ and so that $p< 1+\alpha + {2+\alpha \over \tau^++\epsilon}$ and $\epsilon< {c_1\over 4(-\tau^++N^+)}$,
  since $\varphi$ tends to zero,  there exist positive constant $c, C$ so that for $r_2<r_1$  and    for $r< r_2$ 
                              $$ c r^{-\tau^++\epsilon} \leq u(r) \leq C r^{-\tau^+-\epsilon}$$
                                                        
                              We obtain  with $c_1^+ := {c_1+\epsilon (\tau^+-N^+) \over 1+\epsilon }$, replacing in \eqref{varphi},     
                           $$ ( \varphi (r)r^{-c_1^+})^\prime \leq {C^{p-1-\alpha}\over \Lambda ( \tau^+)^\alpha (1+\epsilon)}   r^{c_2^\prime-{c_1^+  } -1}, $$
                            with $c_2^\prime: = -(p-1-\alpha) (\tau^++\epsilon)+(2+\alpha) >0$.                      We obtain  by integrating ( note that we can choose $\epsilon$ so that $c_1^+\neq c_2^\prime$)
                            $$ \varphi(r_o) r_o^{-c_1^+ }-\varphi(r) r^{-c_1^+ } \leq {C^{p-1-\alpha}\over \Lambda ( \tau^+)^\alpha (1+\epsilon)} {r_o^{c_2^\prime-c_1^+ } -r^{c_2^\prime -c_1^+ }\over c_2^\prime -c_1^+ }\ , $$
                             
                             and then there exists $c_6$ so that 
                             \begin{equation}\label{C5}\varphi(r) \geq -c_6 r^{\inf (c_2^\prime , c_1^+ )}.
                             \end{equation} 
                                                            
                               Let now $\zeta(r) = u(r) r^{\tau^+}$.  We obtain that  for  example if $\varphi(r) \leq 0$ , $\zeta^\prime \leq 0$ and    using \eqref{C5}                         
                                  $$  \zeta^\prime (r) \geq -c_6  r^{\inf (c_2^\prime , c_1^+)-1}\zeta(r) .$$                               Then $\log \zeta$ has a finite limit at zero and then  $\zeta$ has a finite limit $l$ at zero different from zero.   
                                     
                                   If $\varphi \geq 0$ and $\varphi^\prime \geq 0$,   so  are $\zeta$ and $\zeta^\prime$.  In that case we get 
                                   $$\zeta^\prime(r) \leq c_6  r^{\inf (c_2^\prime , c_1^+)-1}\ \zeta.$$
                                   Then $\zeta(0)$ is finite.  In any cases we have 
                                                              $$u(r) \sim  l r^{\tau^+}.$$

               \end{proof}  
                 \begin{proof}of Proposition \ref{phi1}
                
                 Here acting as before ( replacing $\tau^+$ by $\tau^-$ ) one gets 
                 $$(1-(\alpha - \epsilon) { \varphi \over \tau^-})( r\varphi^\prime +\varphi (N^+-2-2\tau^- -\epsilon)+\tau^- ( \tau^--(N^+-2)) )+ {\mu \over \Lambda (\tau^-)^\alpha}  \leq { u^{p-1-\alpha} r^{2+\alpha }\over( \tau^-)^\alpha \Lambda}$$
                   Choosing $r_o$ small enough in order that  for $r< r_o$,
                     $$(1+
|\alpha|) {\varphi(r)\over \tau^-} \leq \epsilon<\inf \left( {1\over 4}, {2+\alpha ( \tau-\tau^-)\over N^+-1-\tau^-}\right),$$
 one derives that 
                     $$ r\varphi^\prime (1-\epsilon) +((2+\alpha ) ( \tau-\tau^-) +\epsilon ( \tau^--N^+)   ) \varphi \leq {u^{p-1-\alpha} r^{2+\alpha}\over( \tau^-)^\alpha \Lambda}$$

 Here since $r\varphi^\prime$ and $\varphi \geq 0$ we can drop $r\varphi^\prime$ and obtain immediately that  for some positive constant $c_3^-$ 
  $$\varphi (r)\leq c_3^-  r^{-(\tau^-+\epsilon) ( p-1-\alpha)+2+\alpha}$$
 
  which is the desired result. We obtain the estimate on $u$ as in the previous proof
  \end{proof}

           \begin{proof} of Theorem\ref{lim=0}
            Let  us recall that  $\beta = {2+\alpha \over p-1-\alpha}$. 
            Let $x(\log r)= r^\beta u(r)$
             and $y(t)  = \beta -{ x^\prime \over x}$.  $y(t) >0$ since $u^\prime(e^t) <0$ as soon as  $t$ is large enough negative. Since $x$ goes to zero,  by Theorem \ref{useconde}, one has $u^{\prime \prime} \geq 0$  around zero and $y$ satisfies \eqref{eqy}.  Using proposition \ref{propy}, we get that $y^\prime$ has a sign, in particular $y^{1+\alpha}$ has a limit. Suppose that 
                    $({y^{1+\alpha}\over 1+\alpha})^\prime >0$, then  $y^{1+\alpha} $ has a finite limit,  and so has  $y$  . Suppose that it is $\neq \tau^\pm$ then there exists $\delta >0$  and $T$ so that  for $t< T$, $\psi$ defined in \eqref{defpsi} satisfies 
                    $\psi( y) > \delta.$ Since $x$ goes to zero, let $T_1< T$ so that  for $t< T_1$, $x^{p-1-\alpha} < { \delta \over 2}$. Then  for $t< T_1$,
                    
                 $$({y^{1+\alpha}\over 1+\alpha})^\prime \geq { \delta \over 2},$$but this contradicts the fact that $y^{1+\alpha} $ has a finite limit. Note that since $({y^{1+\alpha}\over 1+\alpha})^\prime >0$,  the limit cannot be $\tau^-$ since it would imply that $y\in ]\tau^-, \tau^+[$ for $t$ large enough, a contradiction with the equation. We have then   $y\rightarrow \tau^+$. 
                 
                  If $({y^{1+\alpha}\over 1+\alpha})^\prime <0$ suppose first that $y$ tends to $+\infty$.  Then for $t< T$ large enough  negative, 
                  $$ \psi(y) -\Lambda^{-1}x^{p-1-\alpha} \geq { y^{2+\alpha} \over 2}, $$
                   and this contradicts the equation, since $({y^{1+\alpha}\over 1+\alpha})^\prime <0$. Then once more $y^{1+\alpha} $ has a finite limit and  it  cannot be else than $\tau^\pm$. Here the limit can be $\tau^-$.  There remains to prove that $x\sim c e^{(\beta-\tau^\pm)t} $, or that $u(r) \sim c r^{-\tau^\pm}$. 
                   
                   Suppose that we are in the case where $y\rightarrow \tau^+$, 
                  we introduce $\varphi(r) $ so that 
                    $$ {u^\prime r \over u} = -\tau^++\varphi(r)$$
                     Note that if $y\rightarrow \tau^+$ , and $y^\prime \geq 0$, $y\geq \tau^+$, ( fot $t< T$ large enough negative)  and then 
                      both $\varphi$ and $\varphi^\prime $ are negative ( for $r$ small enough ) . If $y^\prime \leq 0$ , $y\leq \tau^+$, and then $\varphi, \varphi^\prime \geq 0$. Moreover $y^\prime$ tends to zero by the equation, which ensures that $r\varphi^\prime \rightarrow 0$.                   Using Proposition \ref{phi},    we get that $u(r) \sim c r^{-\tau^+}$.    
                  
                   In the case $y\rightarrow \tau^-$ we have $y^\prime \leq 0$,     $\varphi$ defined by  $$ {u^\prime r \over u} = -\tau^-+\varphi(r)$$     is  $\geq 0$ as well as $\varphi^\prime$, $r\varphi^\prime$ and $\varphi$ tend to zero,  and then we use Proposition \ref{phi1} and get that 
                   $u(r) \sim c r^{-\tau^-}$.                                          
                   \end{proof}

                    \bigskip
                    
                  We now consider the case where $\limsup u(r) r^\beta >0$, or   the first assertion in item 1:                                     
                                   

                 
                  
 

     \begin{prop}\label{convK}
     
      Suppose that $p< p^\star$ or $p> p^{\star\star}$.       Suppose that $\limsup u(r) r^\beta>0$,  then $\liminf u(r) r^\beta >0$ and finally $u(r) r^\beta $ converges  to $K$
      \end{prop}
      
       \begin{proof}
        Let $l>0$ and $r_n$ so that 
        $u(r_n) r_n^\beta = l_n\rightarrow l$. 
        Let 
        $ \underline{u} = \min ( K, {l\over 2}) r^{-\beta}$.
         It is immediate that $\underline {u}$ is a sub-solution. Furthermore 
         $\underline{u}(r_n) \leq u(r_n)$ and then applying the first comparison principle,  
         $ u \geq \underline{u}$ on $]r_{n+1}, r_n[$. Finally $u \geq \underline{u}$ and $\liminf u >0$. 
         We now prove that $\liminf \geq K$ and $\limsup \leq K$.   Arguing by contradiction,        let $l_1< K$, so that $\liminf u(r) r^\beta = l_1$. 
          Let $\delta$ fixe positif , $\delta < \beta $.            Suppose $\epsilon$ depending on $\alpha, p$ first so that 
           for any $u\in ]0,  2\epsilon[$ one has 
                        \begin{equation}\label{sgnalpha2}1-2| \alpha| u\leq  (1-u)^\alpha \leq 1+2|\alpha|  u, \  
1-2(1+\alpha) |u| \leq (1-u)^{1+\alpha} \leq 1-{1+\alpha\over 2} u
\end{equation} and 

          \begin{equation}\label{alphap}  1-2p u\leq  (1-u)^p \leq 1-{pu\over 2}.
          \end{equation}         Let 
            \begin{equation}\label{cdealpha}
            c( \alpha, \beta , \delta, N^+) = \beta^\alpha ( 1+2 |\alpha|) ( |\beta| | \beta-N^++2|+ | \beta-\delta| | \beta-\delta-N^++2|) + {2\mu(1+\alpha)\over \Lambda}    .          \end{equation} 
            We impose to $\epsilon$ to be also  such that 
            $\epsilon < {1\over 4}$, small in order that \eqref{alphap}, \eqref{sgnalpha2} hold  and 
            $ l_1(1+\epsilon) < { K+ l_1\over 2}$
             and so that 
              \begin{equation}\label{cdealpha}
              2\epsilon ( c(\alpha, \beta , \delta, N^+)-{p\over 2} l_1^{p-1-\alpha}) < K^{p-1-\alpha} -\left({K+l_1\over 2}\right) ^{p-1-\alpha}             \end{equation}
              and $r_o$ so that for $r < r_o$ 
              $ u(r) r^\beta > l_1(1-\epsilon).$
               We take $ a = 2\epsilon r_o^{-\delta}$ and define 
               $$\underline{u} (r) = l_1(1+\epsilon) r^{-\beta} ( 1-a r^\delta).$$
                In particular 
                $$\underline{u} (r_o)   < u(r_o).$$
                By the choice of $\epsilon$ ,  $\underline{u}$ is a sub-solution in $r < r_o$. Indeed, an easy calculation proves that, using \eqref{cdealpha}
                \begin{eqnarray*}
                 | \nabla \underline{u}|^\alpha {\cal M}^+(D^2 u) + \mu  r^{-2-\alpha} \underline{u}^{1+\alpha} &\geq& (  l_1(1+\epsilon) r^{-\beta} )^{1+\alpha} r^{-\alpha-2} ( K^{p-1-\alpha}-  c( \alpha,\beta, \delta,   N^+) a r^\delta) \\
               & \geq& r^{-\beta p} (l_1(1+\epsilon))^p  ( 1-{p\over 2} ar^\delta)\geq   \underline{u}^p
               \end{eqnarray*}
               
                        If $\beta \neq N^+-2$         using the second comparison Theorem \ref{comp2} one gets that 
                 $ u \geq \underline{u} $ which implies in  particular that $\liminf u \geq l_1(1+\epsilon)$,  a contradiction with the hypothesis. 
                 
                  If $\beta = N^+-2$ we observe  that the additional assumption in Theorem \ref{comp2}, say $|u^\prime| \leq c r^{-p-1\over 2+\alpha}$ holds for the super-solution $u$. Indeed, in that case $K^{p-1-\alpha } ={ \mu}$. Then if $l_1< K$
                  $$ { \cal M}^+( D^2 u) \leq -c r^{-\beta p}.$$
                Replacing in radial coordinates,                    $$ ( |u^\prime |^\alpha u^\prime r^{(N^+-1)(1+\alpha)})^\prime \leq -c r^{-1}\leq 0.$$
                In particular $ |u^\prime |^\alpha u^\prime r^{(N^+-1)(1+\alpha)}$ has a limit $\leq 0$ ,  which cannot  be else than finite, since for all $r>0$
             $ |u^\prime |^\alpha (r)u^\prime(r)  r^{(N^+-1)(1+\alpha)}\leq \lim _{r\rightarrow 0} |u^\prime |^\alpha u^\prime r^{(N^+-1)(1+\alpha)}  $.
             We have then $|u^\prime |^\alpha u^\prime r^{(N^+-1)(1+\alpha)}\sim l\leq 0$ which implies that 
              $||u^\prime |^\alpha u^\prime r^{(N^+-1)(1+\alpha)}| \leq C$, 
              or equivalently here 
              $|u^\prime | \leq C^{1\over 1+\alpha} r^{-\beta-1}$.
                                  We end the proof as in the case $\beta \neq  N^+-2$.

                    We now suppose $l_1= \limsup > K$. Let $\epsilon$ so that 
                    $\epsilon < \inf ( {l_1-K\over 2 l_1}, {1\over 4}) $,  and also  so that for $0\leq u\leq 2\epsilon$
             $$ 1+{ | \alpha | \over 2} u\leq (1+u)^\alpha\leq 1+ 2| \alpha| u, \  1+{ 1+\alpha  \over 2} u\leq (1+u)^\alpha \leq 1+ 2(1+\alpha) u,$$
             and 
                 $$  2\epsilon c(\alpha, \beta , \delta, N^+) < \left({K+l_1\over 2}\right) ^{p-1-\alpha}      -K^{p-1-\alpha} $$

                     We take $r_o$ so that for $r < r_o$ 
                     $$u(r) r^\beta \leq l_1 (1+{\epsilon\over 2} ),$$
                      we define 
                      $ a = 2\epsilon  r_o^{-\delta}$ and 
                      $$ \overline{u} =  l_1(1-\epsilon) r^{-\beta} ( 1+ a r^\delta).$$
                       We have 
                       $$\overline{u}(r_o) r_o^\beta= l_1 ( 1-\epsilon) (1+2\epsilon) \geq l_1(1+\epsilon/2)\geq u(r_o) r_o^\beta.$$
                 An analogous calculation as the one done in the case of the sub-solution, proves that $\overline{u}$ is a super-solution in $r< r_o$. Using the second comparison Theorem \ref{comp2}  one gets that $u\leq \overline{u}$ and this  implies $\limsup u (r) r^\beta \leq l_1(1-\epsilon)$ a contradiction.  We have obtained that $\lim u(r) r^\beta=K$. 
              \end{proof}

 \subsection{The case $p \in [1+\alpha+ {2+\alpha \over \tau^+}, 1+\alpha+ {2+\alpha \over \tau^-}[$}
   The proposition below proves item 2 in Theorem    \ref{behav}
     \begin{prop}\label{proplimnulle} Suppose that $p\in [1+\alpha+ {2+\alpha \over \tau^+}, 1+\alpha+ {2+\alpha \over \tau^-}[$. Then 
     $\lim u(r) r^{2+\alpha \over p-1-\alpha} = 0$,  and 
     there exists some positive constant $c$ so that 
     $$u(r) \sim c r^{-\tau^-}$$ 

           \end{prop}
      
       \begin{proof} 
       
        We suppose first that $\beta \neq \tau^+$. Then $\psi(\beta) < 0$.  To prove that $x$ goes to zero, let us suppose by contradiction that $\limsup u(r) r^\beta = l>0$.          Let $0<\delta < \beta$,  Let 
         $$\epsilon <\inf({1\over 4}, c( \alpha, \beta, N^+, \delta ) ^{-1} {| \psi ( \beta)|\over 4})$$ and so that  \eqref{sgnalpha2} and \eqref{alphap} hold
( see the proof of Proposition \ref{convK}), and  $c( \alpha, \beta , N^+,  \delta)$ has been defined in \eqref{cdealpha} .     
         Let  $r_o$ be small enough in order that for $r < r_o$,
               $u(r) r^\beta < l (1+{\epsilon \over 2})$. Let $a = 2\epsilon r_o^{-\delta}$,
         and 
         $\overline{u}(r)=  l(1-\epsilon) r^{-\beta }(1+ a r^\delta)$.
         Then 
         
    $$\overline{u} (r_o) = l   (1-\epsilon)(1+ 2\epsilon)r_o^{-\beta}  > u(r_o).$$
    
     With the assumptions on $\epsilon$, one has  
     $$ | \nabla \overline{u} |^\alpha { \cal M}^+ ( D^2 \overline{u} ) + \mu { \overline{u} ^{1+\alpha} \over r^{2+\alpha} } = r^{-\beta p} ( \psi ( \beta)+ a r^\delta c( \alpha, \beta, N^+, \delta))\leq  { \psi(\beta)\over 2} r^{-\beta p} < 0 < \overline{u} ^p$$
    and then $\overline{u}$ is a super-solution. Then using the second comparison Theorem \ref{comp2} one gets that 
     $u\leq \overline{u}$ in $r< r_o$, and this contradicts $\limsup u(r) r^\beta = l$. We have obtained that $\limsup u r^\beta=  0$.
     
      Using Proposition \ref{propy} one gets that $y^\prime$ has constant sign, hence $y$ converges, and as proved in Theorem \ref{lim=0}, $y\rightarrow \tau^\pm$. We observe that $y$ cannot tend to $\tau^+$. Indeed, if it was the case one would have ${ x^\prime \over x} \sim \beta-\tau^+  <0$, this would imply that $x^\prime \leq -cx$ and integrating one would  get that $x$  would be negative large ,  as soon as  $t$ is  large enough negative, a contradiction. We have then $y\rightarrow \tau^-$ and using 
 proposition \ref{phi1} one gets that $u(r) \sim c r^{-\tau^-}$.                                  
                                                                    
                                                 There remains to do the case $ p = p^\star$, or equivalently $\beta = \tau^+$. 
                                                  In that case $\psi( \beta) = 0$. Suppose that $\limsup u(r) r^{\tau^+} = l>0$. 
                                                  we fix $\delta \in ]0,\beta[ $,    suppose that $\epsilon$ is so that  $\epsilon < {1\over 4}$, 
         $$2\epsilon c( \alpha, \beta, N^+, \delta ) < ({l\over 2})^{p-1-\alpha}, $$
        and  $r_o$ so that for $r< r_o$, 
         $u(r) r^\beta < l (1+{\epsilon \over 2})$. Let $a = 2\epsilon r_o^{-\delta}$, 
         and 
         $\overline{u}(r)=  l(1-\epsilon) r^{-\beta }(1+ a r^\delta)$. 
         Then 
          $$ | \nabla \overline{u} |^\alpha { \cal M}^+ ( D^2 \overline{u} ) + \mu { \overline{u} ^{1+\alpha} \over r^{2+\alpha} }\leq  (l(1-\epsilon))^{1+\alpha} a r^\delta c( \alpha, \beta, N^+) r^{-\beta p} \leq (l(1-\epsilon))^pr^{-\beta p} < \overline{u}^p.$$
        Since 
    $\overline{u} (r_o)  = l   (1-\epsilon)(1+ 2\epsilon) r_o^{-\beta}> u(r_o), $
    using the second comparison Theorem \ref{comp2} one gets also in that case that $\lim u(r) r^\beta = 0$.
     By Proposition \ref{propy} one gets that $y\rightarrow \tau^\pm$. If $y\rightarrow \tau^+$, we cannot have $y^\prime \geq 0$ since this would  imply ${x^\prime \over x} \leq 0$, hence $x^\prime \leq 0$,  and $x$ tends to $0$ and  $x\geq 0$ would imply $x \equiv 0$.      Then $y^\prime \leq 0$, $y< \tau^+$ ( for $t< T$ large enough negative)  and $\varphi > 0$ for $r$ small enough. We have by the proof in proposition \ref{phi}, 
     $$ r\varphi^\prime -c_1^+ \varphi \geq 0$$then $\varphi \leq C r^{c_1^+}$ and this would imply that $u(r) \sim c r^{-\tau^+}$, with $c>0$,  a contradiction with $\lim u(r) r^\beta = 0$. Then $y\rightarrow \tau^-$ and 
using proposition \ref{phi1},  we get 
                                                                                                                $u(r) \sim c r^{-\tau^-}$.

                           \end{proof}

                                  \subsection{The case $ p =     1+ \alpha + {2+\alpha \over \tau^-}$}
                                  
                                   Here we will prove that    
                                 \begin{prop}\label{proptau-}
                                  If $p= 1+ \alpha + {2+\alpha \over \tau^-}$
                                   $$ u(r) \sim  \tilde K r^{-\tau^-} ( -\log r)^{-\tau^-\over 2+\alpha}$$
                                   
                                    with 
                                    $\tilde K = (\Lambda ( \tau^-)^{1+\alpha} ( \tau-\tau^-))^{1\over p-(1+\alpha)}$.
                                    
                                     \end{prop}
                                     
                                     The idea is to compare $u$ with sub- and super-solutions which behave as $r^{-\tau^-} ( -\log r)^{-\tau^-\over 2+\alpha}$ near zero. Note that  here we cannot employ the comparison principle  in Theorem \ref{comp2} since   the bound from below of the  super-solution  is not satisfied. We need another "comparison result " about sub- and super-solutions of an auxiliary equation.

                                     \begin{lemme}\label{lemv}
                                      Suppose that $\beta = p^{\star\star}$ and $u$ is a positive radial solution of \eqref{eqmain}.  
                                     Let $ v = 
                                     r^{\tau^-} u$.
                                    Then ${v^\prime r \over v} \rightarrow 0$, $v^{\prime \prime } r^2 v^{-1}\rightarrow 0$. 
                                     
                                      \end{lemme}                                           \begin{proof} of Lemma \ref{lemv}
                                           
                                       We begin to observe that ${ x^\prime \over x}\rightarrow 0$.                                               For that aim, since $x$ goes to zero, (see proposition \ref{proplimnulle})   we prove as before that $y^\prime $ does not change sign and  then $y$  tends to $\tau^\pm$. 
                                        We prove that it cannot tend to $\tau^+$. If it is the case then 
                                        ${ x^\prime \over x} \sim ( \tau^--\tau^+):= -\gamma <0$.
                                         Then for $t$ large enough negative 
                                                                             one gets the contradiction 
                                          $ x <0$. 
                                          Then $y$ tends to $\tau^-$ and then ${x^\prime \over x}$ tends to zero. 
                                                   Note that this is equivalent to say that ${rv^\prime \over v}$ tends to zero. 
                                                      From the equation \eqref{eqy}, $y^\prime \rightarrow 0$, which means that $r^2 v^{\prime\prime} v^{-1} $
                                                       tends to zero. 
                                                       \end{proof}
                                     
                                         \begin{proof}of Proposition \ref{proptau-} 
                                      Observe that $v(0) = 0$ since $x\rightarrow 0$.
                                                                        An easy calculation proves that  $v$  satisfies  
                                       $$| \tau^- v -rv^\prime |^\alpha ( r^2 v^{ \prime \prime} + ( N^+-1-2\tau^-) v^\prime r + \tau^-  ( \tau^--(N^+-2)) v ) +{\mu\over \Lambda} v^{1+\alpha} ={ v^p\over \Lambda}$$

                                       We prove below that  $v^\prime$ does not change sign, near zero, and since $v(0)=0$ and $v>0$ this  implies  that $v^\prime \geq 0$. 
                                                                             Indeed suppose that we have a sequence $r_n$ coming to zero, so that $v^\prime (r_n) = v^\prime (r_{n+1}) = 0$, and $v^\prime <0$ in $]r_{n+1}, r_n[$,  then by the equation 
                                     $$  ( \tau^- v  )^\alpha ( r^2 v^{ \prime \prime} + \tau^-  ( \tau^--(N^+-2)) v ) +{\mu\over \Lambda} v^{1+\alpha}= ( \tau^-v)^\alpha r^2 v^{\prime \prime} = {v^p\over \Lambda}$$
                                      in particular  $v^{\prime\prime} (r_{n+1})>0$,  a contradiction since $v^{\prime } $ is  negative after  $r_{n+1}$, close to it. 
                                      
                                     {\bf Claim  $v^\prime r v^{-1} \leq Cv^{p-1-\alpha}  $}
                                     
                                      Using ${rv^\prime \over v}$ and ${ r^2 v^{\prime \prime} \over v}$ tend to zero, there exist constants $d_1:=d_1 ( \alpha, N^+, \tau^-)$ , $d_2:= d_2( \alpha, N^+, \tau^-)$, $d_3:=d_3 ( \alpha, N^+, \tau^-)$  $d_4:=d_4 ( \alpha, N^+, \tau^-)$  so that, for any $\epsilon$, there exists  $r_o$ so that $r< r_o$                             
                                                                              $$  ( \tau^- v)^\alpha ( v^{ \prime \prime } r^2+ v^\prime r ( 1+ (2+\alpha) ( \tau-\tau^-)) +  v^\prime r  d_1 \epsilon+ v^{ \prime \prime } r^2d_2 \epsilon ) \leq {v^p\over \Lambda} ,$$
                                                                              and 
                                                                             $$ {v^p\over \Lambda} \leq  ( \tau^- v)^\alpha ( v^{ \prime \prime } r^2+ v^\prime r ( 1+ (2+\alpha) ( \tau-\tau^-)) +  v^\prime r  d_3 \epsilon+ v^{ \prime \prime } r^2d_4 \epsilon ) .$$
                                        
                                         From this one derives assuming $\epsilon$ small enough , but fixe  in order that 
                                         $ \gamma_1 := {1+ d_1 \epsilon + (2+\alpha ) ( \tau-\tau^-)\over 1+ d_2\epsilon } >1$, 
                                         
                                           \begin{equation}\label{gamma2} ( r^{\gamma_1} v^\prime )^\prime \leq {1\over \Lambda ( \tau^-)^\alpha ( 1+ d_2\epsilon)} v^{p-\alpha } r^{\gamma_1-2}
                                           \end{equation}
                                           Using $v^\prime r^{\gamma_1} = v^\prime r r^{ \gamma_1-1} \leq c r^{ \gamma_1-1} v$ one has $v^\prime r^{\gamma_1}\rightarrow 0$ when $r$ tends to zero. 
                                          Integrating, and  using $v$ increasing   one gets 
                                            $$ v^\prime \leq {1\over \Lambda ( \tau^-)^\alpha ( 1+ d_2\epsilon) ( \gamma_1-1)} v^{p-\alpha } r^{-1}.$$

                                               Coming back to \eqref{gamma2}   one also  gets that 
                                               $$|v^{\prime \prime} r^2|\leq c (v^\prime r+ v^{p-\alpha})\leq cv^{p-\alpha}.$$                                                                                                                                                                                To summarize, there exists $d_5$ depending on $N^+, \alpha, \tau^-$ so that, for all $\epsilon$,   there exists some $r_o$ so that   for $r< r_o $ , 
                                                                                      \begin{equation}\label{d3} (1-d_5\epsilon) {v^{p-\alpha} \over (\tau^-)^\alpha \Lambda } \leq v^{\prime \prime } r^2 + (v^\prime r) ( 1+ (2+\alpha)( \tau-\tau^-) ) \leq  (1+d_5\epsilon) {v^{p-\alpha} \over (\tau^-)^\alpha \Lambda } .
                                          \end{equation}

                                           Let $\underline{w} = \underline{K}  ( c_o-\log r)^{-\tau^-\over 2+\alpha}$,          $\underline{K} < \tilde K$ and let $\epsilon $ small enough in order that                             
                                                  $ \underline{K} ^{p-1-\alpha} < { \tilde K^{p-1-\alpha}\over 1+d_5 \epsilon} $, let $r_o$ so that \eqref{d3} holds, and 
                                           $$ c_o = \log r_o+ \left( {  \underline{K} \over v(r_o)}\right)^{2+\alpha \over \tau^-}$$
                                         With that choice, $\underline{w}>0$ in $]0, r_o[$ and   $\underline{w}(r_o) \leq v(r_o)$.
                                            Furthermore  $\underline{w} (0) = 0$ and  it satisfies  for $r< r_o$
                                            $$\underline{w}^{\prime \prime} r^2+ \underline{w}^\prime r(1+ ( N^+-2)(1+\alpha)-(2+\alpha) \tau^-) \geq  \underline{K} \tau^-( \tau-\tau^-) (c_o-\log r)^{{-\tau^-\over 2+\alpha}-1}\geq \underline{w}^{p-\alpha} (1+ d_5\epsilon) \Lambda^{-1}(\tau^-)^{-\alpha}$$
                                           
                                                                                      Let us suppose by contradiction that $\inf (v-\underline{w}) <0$ on $]0, r_o[$. Then let $r_1$ so that the minimum is achieved. We have on $r_1$ , 
                                               $ (v-\underline{w})^\prime (r_1)= 0$, and $ (v-\underline{w})^{\prime\prime}  (r_1)\geq  0$, so
                                               \begin{eqnarray*}
                                              \Lambda^{-1} (1+ d_5\epsilon)  v^{p-\alpha} (r_1)\tau^{-\alpha}&\geq &r_1^2v^{\prime \prime} (r_1)+ r_1v^\prime (r_1)(1+ ( N^+-2)(1+\alpha)-(2+\alpha) \tau^-)  \\
                                                &\geq &  r_1^2 \underline{w}^{\prime \prime} (r_1)+ r_1\underline{w}^\prime (r_1)(1+ (2+\alpha)(\tau- \tau^-)) \\
                                                &\geq&\Lambda^{-1}(1+ d_5\epsilon)  \underline{w}^{p-\alpha} (r_1)(\tau^-)^{-\alpha}
                                                \end{eqnarray*}
                                                a contradiction since $v(r_1) < w(r_1)$.
                                                 We have obtained that 
                                                 $ v\geq \underline{K}  (c_o-\log r)^{-\tau^-\over 2+\alpha}$. 
                                                
                                                To prove some reverse inequality,         let $\overline{K} >  \tilde K$ and let $\epsilon$ small enough in order that    $ \overline{K}^{p-1-\alpha}   >  \tilde K^{p-1-\alpha} (1-2d_5\epsilon)^{-1}$.       Let $r_o$ small in order that \eqref{d3} holds and also so that   ( since $v(0)=0$) \begin{equation}\label{final}
                                               {\tau^-\over 2+\alpha}({\tau^-\over 2+\alpha}+1) \left({v(r_o)\over  \tilde K}\right)^{2+\alpha \over \tau^- } < { d_5\epsilon \over  \Lambda (\tau^-)^\alpha}.
                                               \end{equation}
                                               
                                               We  define                                                 $\overline{w} = \overline{K} (c_o-\log r)^{-\tau^-\over 2+\alpha}$ with 
                                                $$c_o  = \left({ \overline{K} \over v(r_o)}\right)^{2+\alpha\over \tau^-} -(- \log r_o)$$

                                                Then $\overline{w} >0$, is decreasing and convex, $\overline{w}(0)=0$ and $v(r_o) \leq \overline{w}(r_o)$. 
       Furthermore                                                 \begin{eqnarray*}
                                               \overline{w}^{\prime \prime} r^2 &+&  \overline{w}^\prime r(1+ ( N^+-2)(1+\alpha)-(2+\alpha) \tau^-) \\
                                                &=&  \overline{K} \tau^- ( \tau-\tau^-) (c_o-\log r)^{-{\tau^-\over 2+\alpha}-1} 
                                                + \overline{K}   ( {\tau^-\over 2+\alpha}({\tau^-\over 2+\alpha}+1) (c_o-\log r)^{-{\tau^-\over 2+\alpha}-2}\\
                                                                        &\leq &   \overline{K} \tau^- ( \tau-\tau^-) (c_o-\log r)^{-{\tau^-\over 2+\alpha}-1}+ \overline{K}  (c_o-\log r)^{-{\tau^-\over 2+\alpha}-1}  {\tau^-\over 2+\alpha}({\tau^-\over 2+\alpha}+1) ({v(r_o)\over  \tilde K})^{2+\alpha \over \tau^- }  \\
                                                                        &\leq &\overline{K} \tau^- ( \tau-\tau^-) (c_o-\log r)^{-{\tau^-\over 2+\alpha}-1}  (1+ d_5\epsilon)
                                                                                                                                                      \leq (1-d_5\epsilon)  \overline{w}^{p-\alpha} \Lambda^{-1} (\tau^-)^{-\alpha}
                                                \end{eqnarray*}
                                                by using \eqref{final} and the definition of $c_o$. 
                                                                                                 Using the same argument of  comparison  as above but now with $v$ and $\overline{w}$, one obtains that 
                                                  $ v \leq \overline{w}$ and by the arbitrariness of $ \overline{K} $ and $\underline{K}$, one gets the desired equivalent for $v$ and then for $u$.

                                                    \end{proof}
                               \subsection{The case $p> 1+\alpha + {2+\alpha \over \tau^-}$}

                               We  still suppose that $u$ is a positive solution of \eqref{eqmain}, and $x(\log r) = r^\beta u(r)$.
                                                             In  that case we just observe that we cannot have $x\rightarrow 0$. Indeed if it is the case,  $y^\prime $  does not change sign, and then $y$ converges to $\tau^\pm$. But $\beta-\tau^\pm <0$ and this contradicts $y\geq 0$.
                                So $x$ does not tend to zero,  we then use Proposition \ref{convK}. We obtain that $u(r) \sim K r^{-\beta}$ . We could also have used  Theorem of Poincar\'e Bendixson : $x$ must converge to the only stationary  point, which is $K$ ( recall that there is no periodic orbit for $t$ large negative).

                               \end{document}